\documentclass[12pt]{article}

\usepackage[margin=1in]{geometry}
\usepackage{amsmath,amssymb,amsthm}
\usepackage{graphicx}
\usepackage{hyperref}

\usepackage{titlesec}
\titlespacing*{\section}{0pt}{1.2em}{0.8em}
\titlespacing*{\subsection}{0pt}{1em}{0.6em}
\titlespacing*{\subsubsection}{0pt}{0.8em}{0.5em}

\title{\vspace{-1cm}\bfseries A Hierarchy of Geometric Constructions}

\author{MohammadJavad Maarefvand\thanks{Work in progress. This is a preprint version. 
Comments and suggestions are welcome at \href{mailto:javadmaaref80@gmail.com}{javadmaaref80@gmail.com}. 
Licensed under CC BY 4.0 International License. 
Last updated: \today.}}

\begin{document}

\maketitle

\begin{abstract}
This article explores the limits of geometric construction using various tools, both classical and modern. Starting with ruler and compass constructions, we examine how adding methods such as origami, marked rulers (neusis), conic sections, mechanical linkages, and certain transcendental curves expands the range of "constructible" numbers. These methods allow the construction of increasingly complex numbers from square roots, to cube roots, to all algebraic numbers, and in some cases to specific transcendental constants like \(\pi\) and \(e\). We explain how field theory gives a precise way to understand these constructions, and how computability theory shows that no finite geometric method can produce every computable number. In particular, no construction process that can be described step by step can reach uncomputable numbers. The article concludes by presenting a hierarchy of geometric methods, showing how each step increases what is possible, while still leaving strict theoretical limits in place.
\end{abstract}


\tableofcontents
\vspace{1em}
\hrule
\vspace{1em}

\section{Introduction}

For over two thousand years, geometric constructions have raised fundamental questions about the limits of what can be achieved using finite steps and specific tools. Classical geometry, as presented in Euclid’s \emph{Elements}, allowed only the use of a straightedge and compass. By the nineteenth century, the work of Gauss, Wantzel~\cite{Wantzel1837}, and Galois had clearly linked such constructions to towers of quadratic field extensions. As a result, famous problems like angle trisection and cube doubling were shown to be impossible within the strict ruler-and-compass framework.

Later developments expanded the range of allowed geometric tools to include origami folds, marked ruler (neusis) constructions, conic sections, and mechanical linkages. Origami and neusis enable the solution of certain cubic equations, while conic sections can handle some quartic cases. Mechanical linkages, according to Kempe’s Universality Theorem, can in principle produce all algebraic numbers. Going even further, if one assumes access to transcendental curves such as the quadratrix or the logarithmic spiral, it becomes possible to construct specific transcendental constants (such as \(\pi\) and \(e\)) in a finite number of intersection steps.

Despite the increasing power of these methods, recent insights from computability theory show that no finite geometric protocol can construct all computable numbers. Uncomputable numbers remain entirely beyond reach. This article reviews each level of the constructibility hierarchy, from quadratic extensions to constructions involving transcendental curves, and explains why uncountably many real numbers, including many computable ones, fall outside the limits of any finite geometric process.

\subsection{Outline and Main Contributions}

We present an updated hierarchy of geometric constructibility, beginning with classical field theory and extending through origami, conic sections, mechanical linkages, and certain transcendental curves. We also introduce arguments from computability theory to highlight the fundamental limits of finite geometric protocols in capturing all computable real numbers. Section~\ref{sec:RC} reviews ruler-and-compass constructions, followed by origami and neusis in Section~\ref{sec:origami}, conic sections in Section~\ref{sec:conics}, and mechanical linkages in Section~\ref{sec:linkages}. Section~\ref{sec:transcendental} discusses constructions based on transcendental curves. In Section~\ref{sec:additionalcomputability}, we examine computability aspects and explain why even these extended methods remain limited. Section~\ref{sec:modern} briefly covers modern computational approaches. We conclude with a structured overview of all methods, presenting the proposed hierarchy and noting subtle conceptual issues that arise in this context.

\begin{table}[ht!]
\centering
\caption{Examples of Constructible and Non-Constructible Numbers by Method}
\vspace{0.5em}
\label{tab:constructnumbers}
\begin{tabular}{|c|c|c|}
\hline
\textbf{Method} & \textbf{Can Construct} & \textbf{Cannot Construct} \\
\hline
\textbf{Ruler and Compass} 
& \(\sqrt{2}, \tfrac{1}{3}, \cos\!\bigl(\tfrac{2\pi}{17}\bigr)\) 
& \(\sqrt[3]{2}, \cos\!\bigl(\tfrac{2\pi}{9}\bigr), \pi\) \\
\hline
\textbf{Origami / Neusis} 
& \(\sqrt[3]{2}, \cos\!\bigl(\tfrac{2\pi}{9}\bigr), \theta/3\) 
& \(\sqrt[5]{2}, \cos\!\bigl(\tfrac{2\pi}{11}\bigr), \pi\) \\
\hline
\textbf{Conic Sections} 
& Quartic roots (e.g., \(\sqrt[4]{3}\)), \(\cos\!\bigl(\tfrac{2\pi}{15}\bigr)\) 
& \(\sqrt[7]{2}, \cos\!\bigl(\tfrac{2\pi}{23}\bigr), \pi\) \\
\hline
\textbf{Mechanical Linkages} 
& All algebraic numbers (e.g., \(\sqrt[5]{2}, \cos\!\bigl(\tfrac{2\pi}{23}\bigr)\)) 
& \(\pi, e, \ln(2)\) \\
\hline
\textbf{Transcendental Curves} 
& \(\pi, e, \ln(2), \sin(1)\) 
& Chaitin's constant \(\Omega\) \\
\hline
\end{tabular}
\end{table}

\section{Ruler and Compass}
\label{sec:RC}

Ruler-and-compass geometry, as formalized in Euclid’s \emph{Elements}, permits only constructions using a straightedge and compass. In modern algebraic terms, each such construction step corresponds to solving a linear or quadratic equation. As a result, the set of constructible numbers consists precisely of those lying within a tower of quadratic field extensions of \(\mathbb{Q}\). Specifically, a real number \(x\) is constructible by ruler and compass if and only if
\[
\mathbb{Q} = K_0 \subset K_1 \subset \cdots \subset K_n,
\]
where each extension \(K_i\) is obtained by adjoining a square root to \(K_{i-1}\), and \(x \in K_n\). It follows that the degree of the final extension satisfies \([K_n : \mathbb{Q}] = 2^m\) for some integer \(m\).

These results explain why certain classical problems, such as angle trisection and cube doubling, are impossible using only ruler and compass. The solutions to the corresponding cubic equations do not lie within any tower of quadratic extensions. Squaring the circle is ruled out even more strongly, since \(\pi\) is transcendental and cannot occur in any algebraic extension of \(\mathbb{Q}\).

\section{Origami and Neusis}
\label{sec:origami}

Extending ruler-and-compass constructions to include origami (paper folding) or marked ruler (neusis) techniques broadens the range of geometric constructions. Origami, based on the Huzita–Hatori axioms, allows intersections that can solve certain cubic equations in a single fold. Marked ruler constructions, known since the time of Archimedes, similarly enable the extraction of cube roots by sliding a marked straightedge. As a result, both methods make it possible to construct lengths such as \(\sqrt[3]{2}\) and to trisect arbitrary angles. Algebraically, these techniques allow field extensions whose degrees divide \(2^m 3^n\), thereby exceeding the limitations of ruler-and-compass constructions, though still falling short of enabling the solution of general quintic equations and beyond.

\clearpage
\section{Conic Sections}
\label{sec:conics}

Allowing the use of arbitrary conic sections (parabolas, ellipses, and hyperbolas) further increases the power of geometric construction. As early as the 11th century, Omar Khayyám used the intersection of a circle and a parabola to solve certain cubic equations of the form \(x^3 + ax = b\). More complex intersections involving conics can produce solutions to some quartic equations or their special cases, although not all higher-degree polynomials can be solved using a single conic-based construction.

Conic sections therefore allow the solution of a wider range of algebraic problems than origami or neusis alone. However, they still fall short of capturing all algebraic numbers. Their reach typically includes certain quartic equations, but irreducible polynomials of degree five or higher remain unsolvable within this framework.

\section{Mechanical Linkages}
\label{sec:linkages}

Mechanical linkages mark a major step forward in the power of geometric construction. According to Kempe’s Universality Theorem \cite{Kempe1877,Kapovich2002}, any algebraic curve can be traced by a sufficiently intricate arrangement of rigid bars and pivots in the plane. By tracing and intersecting these curves with others (including lines and circles), one can solve any polynomial equation of finite degree, thereby constructing \emph{all} algebraic numbers.

Although physically constructing such linkages for high-degree polynomials is highly complex, the theoretical result remains valid: mechanical linkages can, in principle, generate any algebraic number.

\section{Transcendental Curves}
\label{sec:transcendental}

The final stage in this classical progression assumes that certain transcendental curves, such as the quadratrix, the Archimedean spiral, or the logarithmic spiral, are "given" in the plane. By intersecting these curves with lines or circles in a finite number of steps, one can construct specific transcendental constants.

\subsection*{Constructing \(\pi\) and Other Transcendental Constants}

The quadratrix of Hippias, originally introduced for angle trisection, can also be used to define the constant \(\pi\). Similarly, Archimedean and logarithmic spirals make it possible to construct values such as \(e\) or \(\ln(2)\). The key assumption in these cases is that the entire transcendental curve is pre-defined and available in the plane. This is often justified by a “continuous motion” argument in classical geometry. This assumption diverges from Euclid’s discrete-step model. It effectively allows a continuous process to generate a complete transcendental curve in finite time, after which it may be intersected with elementary curves (such as lines and circles) in a small number of steps.

\clearpage

Even so, allowing a few standard transcendental curves does not make it possible to construct all computable numbers. Many computable but “non-classical” constants cannot be obtained through intersections of a finite family of well-known transcendental curves, and uncomputable reals remain entirely out of reach.

\section{Additional Observations on Computability and Constructibility}
\label{sec:additionalcomputability}

Computability theory offers a deeper understanding of why transcendental curves, despite their expressive power, remain fundamentally limited in what they can construct.

\subsection*{Not All Computable Numbers Become Constructible}

The set of numbers constructible from a \emph{finite} (or even countably infinite but explicitly specified) collection of curves remains far from covering all computable numbers. For example, the positive real solution to the equation
\[
x^x = 2
\]
is computable, as it can be approximated using standard numerical methods. However, there is no guarantee that this number can be obtained by intersecting any finite combination of familiar transcendental curves.

\subsection*{Uncomputable Numbers Cannot Be Constructed}

No finite geometric protocol, defined as a sequence of explicit, algorithmic steps using classical tools (such as ruler, compass, conic intersections, or idealized transcendental curves like the quadratrix, Archimedean spiral, or logarithmic spiral), can produce an uncomputable real number. This is because any such protocol, when described finitely and algorithmically (without non-computable oracles, undecidable properties, or existential choices), can be translated into a Turing machine program via the Church–Turing thesis. Since uncomputable numbers cannot be generated by any algorithmic process, their construction through geometric means would contradict their definition. This applies even to continuous-motion protocols, provided they remain algorithmic and exclude non-computable dependencies.

\subsection*{Hierarchy of Numbers}

The overall \emph{constructibility hierarchy} can be summarized as follows:
\[
\begin{aligned}
\text{(Ruler-Compass Quadratic)} 
&\subsetneq \text{(Origami/Neusis Cubic)} \\
&\subsetneq \text{(Conic Quartic)} \\
&\subsetneq \text{(All Algebraic via Linkages)} \\
&\subsetneq \text{(Some Transcendentals via Transcendental Curves)} \\
&\subsetneq \text{(All Computable Numbers)} \\
&\subsetneq \text{(All Reals).}
\end{aligned}
\]
This progression makes clear that mechanical linkages are sufficient to construct all algebraic numbers, whereas even standard transcendental curves do not suffice to "construct" all computable real numbers. Uncomputable numbers remain strictly beyond the reach of any finite or explicitly defined geometric procedure.

\subsection*{Approximation Versus Exact Construction}

Modern computational tools can approximate the solutions of arbitrary polynomial and transcendental equations to any desired degree of precision. However, such numerical methods are fundamentally different from the classical notion of exact, finite-step geometric construction. As a result, even advanced software does not eliminate the theoretical limitations imposed by field theory, algebra, and computability.

\section{Modern Computational Methods}
\label{sec:modern}

Modern computational tools can solve polynomial and transcendental equations using symbolic and numerical methods. In practice, these methods go well beyond traditional geometric tools by providing high-precision approximations to a wide range of mathematical constants. However, they do not alter the underlying logical framework. Finite-step constructions, even with transcendental curves, remain fundamentally limited, and classical impossibility results still apply.

\section{Conclusion}

Bringing all these methods together yields the following hierarchy of constructibility:

\[
\begin{aligned}[t]
1.\;& \text{Ruler and Compass:} \\
   &\quad \text{Quadratic extensions of } \mathbb{Q} \ (\text{degree } 2^n).\\
2.\;& \text{Origami and Neusis:} \\
   &\quad \text{Field extensions whose degrees divide } 2^m 3^n.\\
3.\;& \text{Conic Sections:} \\
   &\quad \text{Solutions to certain cubic/quartic equations via conic intersections.}\\
4.\;& \text{Mechanical Linkages:} \\
   &\quad \text{All algebraic numbers (Kempe’s Universality).}\\
5.\;& \text{Finite Families of Transcendental Curves:} \\
   &\quad \text{Select transcendental constants (e.g.\ } \pi, e, \ln(2) \text{).}\\
6.\;& \text{All computable numbers:} \\
   &\quad \text{Beyond any single finite or countably definable family of curves.}\\
7.\;& \text{All Real Numbers:} \\
   &\quad \text{Uncountably many, with uncomputables wholly unconstructible.}
\end{aligned}
\]

Ruler-and-compass constructions correspond exclusively to quadratic field extensions, while origami, conic sections, and mechanical linkages enable the construction of algebraic numbers of increasingly higher degrees. Transcendental curves grant access to specific constants such as \(\pi\) and \(e\), but fall far short of capturing the full set of computable real numbers. Uncomputable reals, in turn, remain completely beyond the reach of any finite-step geometric procedure.

\subsection*{Postscript: A Deeper Critique}

\paragraph{(1) Algebraic Foundations.}
The geometric capabilities of ruler-and-compass constructions, along with those of origami and conic methods, can be precisely described within the framework of field theory. Kempe’s original argument regarding mechanical linkages was later formalized in rigorous terms by Kapovich and Millson~\cite{Kapovich2002}, who proved that any algebraic curve can, in principle, be traced by a suitable linkage.

\paragraph{(2) Transcendental Constructions.}
Allowing continuous motions that generate entire transcendental curves in finite time makes certain constants, such as \(\pi\) and \(e\), accessible through finite constructions. However, these are isolated examples. The vast majority of transcendental numbers, including all uncomputable ones and many computable ones, remain beyond reach. Since any finite geometric protocol relies on a fixed set of pre-defined curves, each constructing only countably many values, the total set of numbers constructible in this way is countable. This yields an additional layer in the hierarchy:
\[
\text{Finite Toolbox} \subsetneq \text{Computables} \subsetneq \text{All Reals}.
\]

\paragraph{(3) Uncomputable Numbers.}
Any attempt to construct an uncomputable real is inherently doomed to fail, since any finite specification can, in principle, be transformed into an algorithm, thereby contradicting the very definition of uncomputability. This limitation is absolute.

\paragraph{(4) Defining “Transcendental Geometry.”}
Whereas Euclid’s constructions and mechanical linkages are precisely defined, the notion of “transcendental geometry” is inherently more open-ended. In practice, one selects a finite set of well-known transcendental curves, such as quadratrices or logarithmic and Archimedean spirals, and demonstrates how to construct particular constants from them. However, this toolkit cannot be extended arbitrarily to span all computable or uncomputable numbers.

In the end, the central unifying principle is that finite-step geometric protocols (regardless of their sophistication or the inclusion of continuous curves) remain bound by the fundamental limitations imposed by algebraic field theory and computability.

\clearpage
\section*{References and Further Reading}
\renewcommand{\refname}{}

\end{document}